\documentclass{article}%
\usepackage{amssymb}
\usepackage{amsfonts}
\usepackage{amsmath}
\usepackage{graphicx}%
\providecommand{\U}[1]{\protect\rule{.1in}{.1in}}
\newtheorem{theorem}{Theorem}

\newtheorem{lemma}[theorem]{Lemma}

\newtheorem{proposition}[theorem]{Proposition}

\newenvironment{proof}[1][Proof]{\noindent\textbf{#1.} }{\ \rule{0.5em}{0.5em}}
\begin{document}

\begin{center}
{\large Analytic Continuation of Generalized Trigonometric Functions}

\bigskip

Pisheng Ding

Illinois State University

Normal, Illinois, 61790

\texttt{pding@ilstu.edu}

\bigskip
\end{center}

\begin{quote}
\textsc{Abstract}. Via a unified geometric approach, a class of generalized
trigonometric functions with two parameters are analytically extended to
maximal domains on which they are univalent. Some consequences are deduced
concerning radius of convergence for the Maclaurin series, commutation with
rotation, continuation beyond the domain of univalence, and periodicity.
\end{quote}



\section{Introduction}

Inverses of functions of the form $y\mapsto\int_{0}^{y}(1-x^{q})^{-1/p}dx$ for
$y\in\lbrack0,1]$ with $p,q>1$ have been of interest to analysts. See, for
example, \cite{Bh}, \cite{Bind}, \cite{Ed}, \cite{L-E}; also see \cite{L} for
an account of early work in this area. Herein, we study the complex-analytic
aspects of a subclass of them and identify their maximal domain of univalence.

Throughout this note, $n$ and $k$ are integers with $n>2$ and $1\leq k<n$.

For $y\in\lbrack0,1]$, let%
\[
F_{n,k}(y)=\int_{0}^{y}\frac{1}{(1-x^{n})^{k/n}}dx
\]
Denoting the number $F_{n,k}(1)$ by $\varphi_{n,k}$, define $S_{n,k}%
:[0,\varphi_{n,k}]\rightarrow\lbrack0,1]$ to be the inverse of $F_{n,k}$.
These functions $S_{n,k}$ are often referred to as generalized sine functions.
Primarily studied as real-valued functions, they are sometimes considered as
analytic functions on $F_{n,k}[\mathbb{D}]$ ($\mathbb{D}=\{z\in%
\mathbb{C}
\mid|z|<1\}$) due to the fact that $F_{n,k}$ is univalent as a complex-valued
function on $\mathbb{D}$ (by the Noshiro-Warschawski theorem in \cite{Duren}).

Note that, in the literature, the function $S_{n,k}$ is more commonly denoted
by $\sin_{n/k,n}$, whereas the number $2\varphi_{n,k}$ is often notated as
$\pi_{n/k,n}$.

In this note, we identify, for each $n$ and $k$, $S_{n,k}$'s \textquotedblleft
natural\textquotedblright\ domain of analyticity, which turns out to be the
maximal domain on which $S_{n,k}$ is univalent. Our unified treatment also
encompasses the circular sine function (if we allow $n=2$), the historically
important lemniscate sine function $S_{4,2}$, and the Dixon's elliptic
function $S_{3,2}$.

In \S 2, we introduce notation and state the main results, which are then
proved in \S 3. In \S 4, we note several consequences, some of which concern
further analytic continuation on larger domains beyond the domain of univalence.

\section{The Main Results}

Before we state a key lemma and the main theorem, we introduce notation.

We often denote points in $%
\mathbb{C}
$ by capital letters in the English alphabet; when we do so, we write $AB$ for
the (closed) line segment between point $A$ and point $B$, whereas we write
$[A,B)$ for the half-open segment $AB\setminus\{B\}$.

Now, let $n>2$ be fixed. Let $\omega_{n}=e^{i2\pi/n}$. Denote the point
$\varphi_{n,k}$ by $A_{n,k}$ and the point $\varphi_{n,k}\omega_{n}$ by
$B_{n,k}$. We construct a closed set $\Pi_{n,k}$ for each $k<n$. There are two
cases, $k=1$ and $k>1$, that need separate (but related) treatments.

We first treat the case $k=1$ and construct $\Pi_{n,1}$. The boundary of
$\Pi_{n,1}$ is the union of the two rays $\left\{  \varphi_{n,1}+te^{i\pi
/n}\mid t\geq0\right\}  $ and $\left\{  \varphi_{n,1}\omega_{n}+te^{i\pi
/n}\mid t\geq0\right\}  $ along with the two segments $OA_{n,1}$ and
$OB_{n,1}$, whereas the interior of $\Pi_{n,1}$ is the component of $%
\mathbb{C}
\smallsetminus\partial\Pi_{n,1}$ that contains $\{te^{i\pi/n}\mid t>0\}$ (the
bisector of $\angle A_{n,1}OB_{n,1}$.)

For $k>1$, let $P_{n,k}$ denote the point on the bisector of $\angle
A_{n,k}OB_{n,k}$ such that $A_{n,k}P_{n,k}$ has an angle of inclination
$k\pi/n$. (If $k$ were $1$, such a point would only \textquotedblleft
exist\textquotedblright\ at $\infty$; cf. the construction of $\Pi_{n,1}$.) We
denote by $\Pi_{n,k}$ the (compact) set enclosed by the polygon $OA_{n,k}%
P_{n,k}B_{n,k}$. Note that $\measuredangle OP_{n,k}A_{n,k}=(k-1)\pi/n$. Thus,
in $\Pi_{n,k}$, the interior angle at the vertex $P_{n,k}$ exceeds $\pi$ iff
$k>\frac{n}{2}+1$, in which case $\Pi_{n,k}$ is nonconvex, whereas the
\textquotedblleft interior angle\textquotedblright\ at $P_{n,k}$ becomes a
straight angle (and $P_{n,k}$ is a degenerate vertex) iff $k=\frac{n}{2}+1$.
With a little geometry, one can verify that $P_{n,k}$ is the point%
\begin{equation}
\varphi_{n,k}e^{i\pi/n}\left(  \cos\frac{\pi}{n}+\sin\frac{\pi}{n}\cot
\frac{(k-1)\pi}{n}\right)  \text{ ;} \label{P}%
\end{equation}
note that $\cot\left[  (k-1)\pi/n\right]  <0$ when $k>\frac{n}{2}+1$, i.e.,
when the interior angle at the vertex $P_{n,k}$ exceeds $\pi$.

Let $V_{n}=\{re^{i\theta}\mid r>0;\,\theta\in(0,2\pi/n)\}$. We state a key lemma.

\begin{lemma}
\label{Key}The analytic continuation of $F_{n,k}$ is a conformal equivalence
from $V_{n}$ to $\mathring{\Pi}_{n,k}$ (the interior of $\Pi_{n,k}$); its
continuous extension to $\overline{V}_{n}$ (the closure of $V_{n}$) restricts
to a homeomorphism from $\partial V_{n}$ onto $\partial\Pi_{n,k}$ when $k=1$
and onto $\partial\Pi_{n,k}\setminus\{P_{n,k}\}$ when $k>1$.
\end{lemma}

To state the main Theorem \ref{Main}, define the domain $\Omega_{n,k}$ to be
the interior of%
\[
\bigcup_{j=0}^{n-1}\left(  \omega_{n}^{j}\cdot\Pi_{n,k}\right)  \text{ ;}%
\]
letting $J=[1,\infty)$, define $\Sigma_{n}$ to be the plane with the $n$ slits
$\omega_{n}^{j}\cdot J$ ($j\in%
\mathbb{Z}
$):%
\[
\Sigma_{n}=%
\mathbb{C}
\setminus\left(  \bigcup_{j=0}^{n-1}\left(  \omega_{n}^{j}\cdot J\right)
\right)  \text{.}%
\]

\begin{theorem}
\label{Main}Concerning $S_{n,k}$, we have the following statements.

\begin{enumerate}
\item The analytic continuation of $S_{n,k}$ is a conformal equivalence from
$\Omega_{n,k}$ onto $\Sigma_{n}$.

\item $S_{n,k}:\Omega_{n,k}\rightarrow\Sigma_{n}$ has a continuous extension
$\tilde{S}_{n,k}$ with $\tilde{S}_{n,k}(\omega_{n}^{j}\varphi_{n,k}%
)=\omega_{n}^{j}$ ($j\in%
\mathbb{Z}
$); moreover,

\begin{description}
\item[(a)] $\tilde{S}_{n,1}$ maps $\partial\Omega_{n,1}\setminus\cup
_{j}\{\omega_{n}^{j}\varphi_{n,1}\}$ two-to-one onto $\cup_{j}\omega_{n}%
^{j}\cdot\mathring{J}$;

\item[(b)] for $k>1$, $\tilde{S}_{n,k}$ maps $\partial\Omega_{n,k}%
\setminus\cup_{j}\{\omega_{n}^{j}\varphi_{n,k}\,,\,\omega_{n}^{j}\cdot
P_{n,k}\}$ two-to-one onto $\cup_{j}\omega_{n}^{j}\cdot\mathring{J}$.
\end{description}

\item $\Omega_{n,k}$ is the maximal domain containing $0$ on which $S_{n,k}$
is univalent.
\end{enumerate}
\end{theorem}

We now set out to make the case.

\section{Analytic Continuation of $F$ and $S$; Proof of Main Results}

Let%
\[
K_{n,k}(z)=\frac{1}{(1-z^{n})^{k/n}}%
\]
with the requirement that $K_{n,k}(0)=1$ and that $K_{n,k}$ be continuous on
$\overline{V}_{n}\setminus\{1,\omega_{n}\}$. Then, $K_{n,k}$ is analytic on
$V_{n}$ with a primitive%
\[
F_{n,k}(z)=\int_{0}^{z}K_{n,k}(\zeta)d\zeta\text{ ,}%
\]
where the integral is path-independent. We first examine the behavior of
$F_{n,k}$ on $\partial V_{n}$.

In the following, $\sqrt[n]{a}$ stands for the principal $n$th root of a
\textit{nonnegative} number $a$. Also recall that $J:=[1,\infty)$.

\begin{itemize}
\item For $x\in\lbrack0,1)$, $F_{n,k}(x)=\int_{0}^{x}1/\sqrt[n]{1-t^{n}}%
^{k}dt$. Thus,
\[
F_{n,k}[[0,1]]=OA_{n,k}\text{ .}
\]

\item For $x>1$, $K_{n,k}(x)=e^{ik\pi/n}/\sqrt[n]{x^{n}-1}^{k}$ where the
phase factor $e^{ik\pi/n}$ is dictated by the continuity of $K_{n,k}$. Then,%
\[
F_{n,k}(x)=\varphi_{n,k}+e^{ik\pi/n}\int_{1}^{x}\frac{1}{\sqrt[n]{t^{n}-1}%
^{k}}dt\text{.}%
\]

In the case $k=1$, since $\lim_{x\rightarrow\infty}\int_{1}^{x}\left(
\sqrt[n]{t^{n}-1}\right)  ^{-1}dt=\infty$,%
\[
F_{n,1}[J]=\left\{  \varphi_{n,1}+te^{i\pi/n}\mid t\geq0\right\}  \text{ ,}%
\]
which is a ray originating from $A_{n,1}$ with an angle of inclination $\pi/n$.

If $k>1$, then%
\[
\lim_{x\rightarrow\infty}F_{n,k}(x)=\varphi_{n,k}+e^{ik\pi/n}I_{n,k}%
\text{\quad where }I_{n,k}=\int_{1}^{\infty}\frac{1}{\sqrt[n]{t^{n}-1}^{k}%
}dt\text{ .}%
\]
and%
\[
F_{n,k}[J]=\left\{  \varphi_{n,k}+te^{ik\pi/n}\mid t\in\lbrack0,I_{n,k}%
)\right\}  \text{ ,}%
\]
which is a segment $[A_{n,k},P_{n,k}^{\prime})$ of length $I_{n,k}$ with angle
of inclination $k\pi/n$.

\item For $z=t\omega_{n}$ with $t\in\lbrack0,1)$,%
\[
K_{n,k}(z)=(1-(t\omega_{n})^{n})^{-k/n}=\left(  \sqrt[n]{1-t^{n}}\right)
^{-k}%
\]
and%
\[
F_{n,k}(z)=\int_{0}^{t}K_{n,k}(\tau\omega_{n})\omega_{n}d\tau=\omega_{n}%
\int_{0}^{t}\frac{1}{\sqrt[n]{1-\tau^{n}}^{k}}d\tau\text{.}%
\]
Thus,%
\[
F_{n,k}[\{t\omega_{n}\mid t\in\lbrack0,1]\}]=OB_{n,k}\text{ .}%
\]

\item For $z=t\omega_{n}$ with $t>1$,%
\[
K_{n,k}(z)=e^{-i\pi k/n}\left(  \sqrt[n]{t^{n}-1}\right)  ^{-k}%
\]
where the phase factor is again dictated by the continuity of $K_{n,k}$, and%
\begin{align*}
F_{n,k}(z)  &  =F_{n,k}(\omega_{n})+\int_{1}^{t}K_{n,k}(\tau\omega_{n}%
)\omega_{n}d\tau\\
&  =\varphi_{n,k}\omega_{n}+e^{-ik\pi/n}\omega_{n}\int_{1}^{t}\frac
{1}{\sqrt[n]{\tau^{n}-1}^{k}}d\tau\text{ .}%
\end{align*}

In the case $k=1$,%
\[
F_{n,1}[\omega_{n}\cdot J]=\left\{  \varphi_{n,1}\omega_{n}+te^{i\pi/n}\mid
t\geq0\right\}
\]
which is a ray originating from $B_{n,1}$ with an angle of inclination $\pi/n$.

If $k>1$, then,%
\[
F_{n,k}[\omega_{n}\cdot J]=\left\{  \varphi_{n,k}\omega_{n}+t\omega
_{n}e^{-ik\pi/n}\mid t\in\lbrack0,I_{n,k})\right\}
\]
is a segment $[B_{n,k},P_{n,k}^{\prime\prime})$ of length $I_{n,k}$ whose
angle with $\overrightarrow{OB}_{n,k}$, measured \textit{clockwise from
}$\overrightarrow{OB}_{n,k}$, has measure $k\pi/n$.

\item We claim that, for $k>1$,%
\begin{equation}
\lim_{z\in\overline{V}_{n},\;|z|\rightarrow\infty}F_{n,k}(z)=\lim
_{r\rightarrow\infty}F_{n,k}(r)\text{ ,} \label{F(infinity)}%
\end{equation}
which will imply that $P_{n}^{\prime\prime}=P_{n}^{\prime}=P_{n}$. Let
$\epsilon>0$ be given. Note that%
\[
F_{n,k}(re^{i\theta})-F_{n,k}(r)=\int_{\zeta=re^{it},\;t\in\left[
0,\theta\right]  }\frac{1}{(1-\zeta^{n})^{k/n}}d\zeta\text{ ,}%
\]
whose modulus can be bounded by $\epsilon/2$ for sufficiently large $r$ (since
$k>1$) and for all $\theta\in\left[  0,2\pi/n\right]  $. Also $|F_{n,k}%
(r)-P_{n,k}^{\prime}|<\epsilon/2$ for sufficiently large $r$ as $P_{n,k}%
^{\prime}=\lim_{r\rightarrow\infty}F_{n,k}(r)$. Thus, $|F_{n,k}(re^{i\theta
})-P_{n,k}^{\prime}|<\epsilon$ for all sufficiently large $r$ and for all
$\theta\in\left[  0,2\pi/n\right]  $, proving the claim.
\end{itemize}

\medskip

\noindent\textbf{Remark}. Let $k>1$. Consider $\triangle OA_{n,k}P_{n,k}$. On
one hand, the point $P_{n,k}$ is located by (\ref{P}); on the other hand,
since $|A_{n,k}P_{n,k}|=I_{n,k}$, $P_{n,k}$ is also the point%
\[
\left(  \varphi_{n,k}\cos\frac{\pi}{n}+I_{n,k}\cos\frac{(k-1)\pi}{n}\right)
e^{i\pi/n}\text{ .}%
\]
Comparing the two, we obtain the identity%
\[
\int_{1}^{\infty}\frac{1}{(t^{n}-1)^{k/n}}dt=\frac{\sin(\pi/n)}{\sin
[(k-1)\pi/n]}\varphi_{n,k}\,\text{.}%
\]

\medskip

Returning to the analysis of $F_{n,k}$ on $V_{n}$, we claim that $F_{n,k}$
maps $V_{n}$ bijectively onto $\mathring{\Pi}_{n,k}$. Let $w\in\mathring{\Pi
}_{n,k}$. Let $\gamma_{r}$ be the positively oriented boundary of the circular
sector with vertices $0$, $r$, and $r\omega_{n}$. The preceding analysis shows
that, for all sufficiently large $r$, $F_{n,k}\circ\gamma_{r}$ winds around
$w$ exactly once.

Thus, $F_{n,k}:V_{n}\rightarrow\mathring{\Pi}_{n,k}$ is a conformal
equivalence. The boundary behavior of $F_{n,k}$ detailed above shows that its
continuous extension is a homeomorphism from $\partial V_{n}$ onto
$\partial\Pi_{n,k}$ when $k=1$ and onto $\partial\Pi_{n,k}\setminus
\{P_{n,k}\}$ when $k>1$.

This concludes the argument for Lemma \ref{Key}.

At long last, we define $S_{n,k}:\mathring{\Pi}_{n,k}\rightarrow V_{n}$ to be
the inverse of $F_{n,k}$. The boundary extension of $S_{n,k}$ mirrors that of
$F_{n,k}$ and maps $OA_{n,k}$ and $OB_{n,k}$ to the two line segments $[0,1]$
and $\{t\omega_{n}\mid t\in\lbrack0,1]\}$ on $\partial V_{n}$. We may then
apply the Schwarz reflection principle repeatedly to analytically continue
$S_{n,k}$ over $\Omega_{n,k}$ with range $\Sigma_{n}$. This establishes Part 1
of Theorem \ref{Main}

Note that $S_{n,k}[\omega_{n}^{j}\cdot\mathring{\Pi}_{n,k}]=\omega_{n}%
^{j}\cdot V_{n}$. By the boundary behavior of $F_{n,k}$, the boundary
extension of $S_{n,k}$ behaves as described by Part 2 of Theorem \ref{Main}.

Finally, we argue Part 3 of Theorem \ref{Main}, i.e., that $\Omega_{n,k}$ is
maximal among domains on which $S_{n,k}$ is univalent. Suppose that $S_{n,k}$
is analytically continued on an (open connected)\ domain properly containing
$\Omega_{n,k}$. This domain necessarily contains a disc $U$ around some line
segment $L\subset\partial\Omega_{n,k}$ with $L$ mapped by $\tilde{S}_{n,k}$
into a ray $\omega_{n}^{j}\cdot J$. By the Schwarz reflection principle,
$S_{n,k}$ maps $U\setminus(\Omega_{n,k})$ into $\Sigma_{n}=S_{n,k}%
(\Omega_{n,k})$, ruining univalence.

\section{Some Consequences}

From the main results, we deduce some notable consequences.

First, we consider the radius of convergence of the Maclaurin series for
$S_{n,k}$.

\begin{proposition}
\label{Radius}Let $R_{n,k}$ be the radius of convergence for the Maclaurin
series for $S_{n,k}$. Then,

\begin{enumerate}
\item $R_{n,1}=\varphi_{n,1}$;

\item for $k>1$, $R_{n,k}\leq\left(  \cos\frac{\pi}{n}+\sin\frac{\pi}{n}%
\cot\frac{(k-1)\pi}{n}\right)  \varphi_{n,k}$;

\item for $k\geq\frac{n}{2}+1$, $R_{n,k}=\left(  \cos\frac{\pi}{n}+\sin
\frac{\pi}{n}\cot\frac{(k-1)\pi}{n}\right)  \varphi_{n,k}$.
\end{enumerate}
\end{proposition}

\begin{proof}
Recall that, by (\ref{P}), $|OP_{n,k}|=\left(  \cos\frac{\pi}{n}+\sin\frac
{\pi}{n}\cot\frac{(k-1)\pi}{n}\right)  \varphi_{n,k}$.

\begin{enumerate}
\item Note that $\operatorname*{dist}(0,\,%
\mathbb{C}
\setminus\Omega_{n,1})=\varphi_{n,1}$ and therefore $R_{n,1}\geq\varphi_{n,1}%
$. We claim that no extension of $S_{n,1}$ can be analytic at $\varphi_{n,1}$.
To see this, first note that the two rays $\{\varphi_{n,1}+te^{i\pi/n}\mid
t\geq0\}$ and $\{\varphi_{n,1}+te^{-i\pi/n}\mid t\geq0\}$ are mapped by
$\tilde{S}_{n,1}$ onto $J$ with $\tilde{S}_{n,1}(\varphi_{n,1})=1$. Within
$\Omega_{n,1}$, these two rays make an interior angle of measure $2\pi
(n-1)/n$. Thus, $\tilde{S}_{n,1}$ expands angle at $\varphi_{n,1}$ by the
non-integer factor $n/(n-1)$ and hence cannot be analytic there. (If we allow
$n=2$, $n/(n-1)$ is an integer and indeed the ordinary circular sine function
\textit{is} analytic at $\pi/2$!)

\item Let $k>1$. Since $|S_{n,k}(z)|\rightarrow\infty$ as $z\rightarrow
P_{n,k}$ in $\Omega_{n,k}$, $R_{n,k}\leq|OP_{n,k}|$.

\item When $k\geq\frac{n}{2}+1$, the interior angle of the polygon
$\overline{\Omega}_{n,k}$ at $P_{n,k}$ is at least $\pi$ and hence
$|OP_{n,k}|=\operatorname*{dist}(0,\,%
\mathbb{C}
\setminus\Omega_{n,k})$. Hence, $R_{n,k}\geq|OP_{n,k}|$.
\end{enumerate}
\end{proof}

Next we examine the interaction between $S_{n,k}$ and certain rotations around
$O$. We need a preliminary observation. Let $L_{n}=\{te^{i\pi/n}\mid t\geq0\}$.

\begin{lemma}
\label{Lemma L}When $k>1$, $F_{n,k}[L_{n}]=[O,P_{n,k})\subset L_{n}$, whereas
$F_{n,1}[L_{n}]=L_{n}$.
\end{lemma}

\begin{proof}
Integrating $K_{n,k}$ along $L_{n}$, we obtain%
\[
F_{n,k}(te^{i\pi/n})=e^{i\pi/n}\int_{0}^{t}\frac{1}{\sqrt[n]{1+\tau^{n}}^{k}%
}d\tau\text{ .}%
\]
The claim for $k>1$ follows at once in light of (\ref{F(infinity)}), whereas
the claim about $F_{n,1}$ is due to the divergence of $\int_{0}^{\infty
}1/\sqrt[n]{1+\tau^{n}}\,d\tau$.
\end{proof}

\medskip

\noindent\textbf{Remarks.} Let $k>1$.

\begin{itemize}
\item The consideration given for Lemma \ref{Lemma L} yields another
expression for $|OP_{n,k}|$, i.e., $\int_{0}^{\infty}\left(  1+t^{n}\right)
^{-k/n}dt$, in addition to that given by (\ref{P}). Comparing the two, we
obtain the identity%
\[
\int_{0}^{\infty}\frac{1}{\left(  1+t^{n}\right)  ^{k/n}}dt=\left(  \cos
\frac{\pi}{n}+\sin\frac{\pi}{n}\cot\frac{(k-1)\pi}{n}\right)  \varphi
_{n,k}\text{ .}%
\]

\item Define $\sqrt{V_{n}}$ to be $\{re^{i\theta}:r>0;\theta\in(0,\pi/n)\}$.
By Lemma \ref{Lemma L}, $F_{n,k}$ maps $\sqrt{V_{n}}$ conformally onto
$\mathring{\triangle}OA_{n,k}P_{n,k}$, the interior of $\triangle
OA_{n,k}P_{n,k}$. There is a \textit{unique} conformal equivalence%
\[
\Psi_{n,k}:\{z\mid\operatorname{Im}z>0\}\rightarrow\mathring{\triangle
}OA_{n,k}P_{n,k}%
\]
whose continuous extension maps $%
\mathbb{R}
\cup\{\pm\infty\}$ onto $\partial(\triangle OA_{n,k}P_{n,k})$ with $\Psi
_{n,k}(0)=O$, $\Psi_{n,k}(1)=A_{n,k}$, and $\Psi_{n,k}(\pm\infty)=P_{n,k}$. By
uniqueness, $F_{n,k}(z)=\Psi_{n,k}(z^{n})$ for $z\in\sqrt{V_{n}}$. As
$\Psi_{n,k}$ can be expressed by a Schwarz-Christoffel integral formula, this
functional identity yields another integral expression for $F_{n,k}$ on
$\sqrt{V_{n}}$.
\end{itemize}

\medskip

We now show that each $S_{n,k}$ commutes with rotation around $O$ by angle
$2\pi/n$.

\begin{proposition}
\label{Prop Rotation}For $z\in\Omega_{n,k}$, $S_{n,k}(\omega_{n}z)=\omega
_{n}S_{n,k}(z)$.
\end{proposition}

\begin{proof}
For a line $L$ in the plane, let $R_{L}$ denote the reflection across $L$.
Recall that the composition of reflections across two intersecting lines is a
rotation around their point of intersection by an angle that is twice the
angle between the two lines.

It suffices to check this identity for $z$ in the interior of $\triangle
OA_{n,k}P_{n,k}$. By Lemma \ref{Lemma L}, $S_{n,k}$ maps $[O,P_{n,k})\subset
L_{n}$ into $L_{n}$. Applying the Schwarz reflection principle, we obtain%
\[
S_{n,k}\left(  R_{L_{n}}(z)\right)  =R_{L_{n}}\left(  S_{n,k}(z)\right)
\text{.}%
\]
Applying this and the definition of $S_{n,k}$ on $\omega_{n}\cdot\mathring
{\Pi}_{n,k}$, we have%
\begin{align*}
\left(  S_{n,k}\circ R_{(\omega_{n}\cdot%
\mathbb{R}
)}\right)  (R_{L_{n}}(z))  &  =\left(  R_{(\omega_{n}\cdot%
\mathbb{R}
)}\circ S_{n,k}\right)  (R_{L_{n}}(z))\\
&  =\left(  R_{(\omega_{n}\cdot%
\mathbb{R}
)}\circ R_{L_{n}}\right)  \left(  S_{n,k}(z)\right)  \text{ ;}%
\end{align*}
i.e.,%
\[
S_{n,k}\circ\left(  R_{(\omega_{n}\cdot%
\mathbb{R}
)}\circ R_{L_{n}}\right)  =\left(  R_{(\omega_{n}\cdot%
\mathbb{R}
)}\circ R_{L_{n}}\right)  \circ S_{n,k}\text{ .}%
\]
Because $R_{\omega_{n}\cdot%
\mathbb{R}
}\circ R_{L_{n}}(\zeta)=\omega_{n}\cdot\zeta$ for all $\zeta\in%
\mathbb{C}
$, the result follows.
\end{proof}

\medskip

Next we consider the possibility of further continuation of $S_{n,k}$.

\begin{proposition}
\label{Prop: n Even}Suppose that $n=2k$. Then $S_{n,k}$ can be analytically
continued to a function (also denoted by $S_{n,k}$) on the interior of%
\[
\cup_{j=0}^{n-1}\left[  \omega_{n}^{j}\cdot\left(  \cup_{m\in%
\mathbb{Z}
}\left(  2m\varphi_{n,k}+\overline{\Omega}_{n,k}\right)  \right)  \right]
\text{ ;}%
\]
for $z\in\omega_{n}^{j}\cdot\left(  \cup_{m\in%
\mathbb{Z}
}\left(  2m\varphi_{n,k}+\overline{\Omega}_{n,k}\right)  \right)  $,%
\[
S_{n,k}(z+4\omega_{n}^{j}\varphi_{n,k})=S_{n,k}(z)\text{ .}%
\]
Furthermore, $S_{n,k}$ has multiplicity 2 at $\omega_{n}^{j}\varphi_{n,k}$ for
$j\in\{0,1,\cdots,(n-1)\}$.
\end{proposition}

\begin{proof}
When $k=n/2$, $\measuredangle OA_{n,k}P_{n,k}=\pi/2$. Therefore, in the
polygon $\overline{\Omega}_{n,k}$, the interior angle at $A_{n,k}$ becomes a
straight angle. Let $P_{n,k}^{\ast}$ denote the complex conjugate of $P_{n,k}%
$. Note that $\tilde{S}_{n,k}$ folds up the open segment $(P_{n,k}%
,P_{n,k}^{\ast})$ into $J$ with $\tilde{S}_{n,k}(A_{n,k})=1$. We may apply the
Schwarz reflection principle to continue $S_{n,k}$ across $(P_{n,k}%
,P_{n,k}^{\ast})$, and then across $2\varphi_{n,k}+(P_{n,k},P_{n,k}^{\ast})$,
and so on. Recall that the composition of reflections across two parallel
lines is a translation by twice the distance between them. This implies that
$4\varphi_{n,k}$ is a period for $S_{n,k}$ on $\cup_{m\in%
\mathbb{Z}
}\left(  2m\varphi_{n,k}+\overline{\Omega}_{n,k}\right)  $. Applying the
identity in Proposition \ref{Prop Rotation}, we can extend $S_{n,k}$
analytically on the domain claimed by the statement. Finally, note that
$S_{n,k}$ is two-to-one on some disc centered at $A_{n,k}$, proving the final claim.
\end{proof}

\medskip

\noindent\textbf{Remark}. Note that $S_{4,2}$ is the historically important
lemniscate sine function $\operatorname*{sl}$ (and $2\varphi_{4,2}$ is the
lemniscate constant). Since $\overline{\Omega}_{4,2}$ is a square each of
whose sides is mapped into $%
\mathbb{R}
$ or $i%
\mathbb{R}
$, repeated application of Schwarz reflection principle allows $S_{4,2}$ to be
analytically continued to an elliptic function on $%
\mathbb{C}
$, a well-known classical result now encompassed by Proposition
\ref{Prop: n Even}.

\medskip

By considering the behavior of $\tilde{S}_{n,k}$ near $P_{n,k}$, we deduce the following.

\begin{proposition}
\label{(n/2)th Power}Suppose that $n$ is even and $k=\frac{n}{2}+1$. Then,
$(S_{n,k})^{n/2}$ can be continued to a meromorphic function on a neighborhood
of $P_{n,k}$, whereas $S_{n,k}$ cannot.
\end{proposition}

\begin{proof}
When $k=\frac{n}{2}+1$, $P_{n,k}$ becomes a degenerate vertex of
$\overline{\Omega}_{n,k}$ as the interior angle at $P_{n,k}$ becomes a
straight angle. Recall that $\tilde{S}_{n,k}$ maps $[A_{n,k},P_{n,k})$ onto
$J$ and that it maps $[B_{n,k},P_{n,k})$ onto the ray $\omega_{n}\cdot J$. If
$S_{n,k}$ were continued to a meromorphic function on a disc around $P_{n,k}$,
then Schwarz reflection across $[A_{n,k},P_{n,k})$ and across $[B_{n,k}%
,P_{n,k})$ would entail a contradiction.

Note that $(S_{n,k})^{n/2}$ maps $[B_{n,k},P_{n,k})$ onto $(-\infty,1]$. Thus,
$[A_{n,k},B_{n,k}]\setminus\{P_{n,k}\}$ is mapped into the $x$-axis by
$(S_{n,k})^{n/2}$. The Schwarz reflection principle can then be applied to
analytically extend $(S_{n,k})^{n/2}$ on some punctured disc centered at
$P_{n,k}$. Thus, $P_{n,k}$ is an isolated singularity of $(S_{n,k})^{n/2}$ and
in fact a pole.
\end{proof}

\section{Postscript}

\noindent\textbf{Remarks.}

\begin{enumerate}
\item Our main results unify various individual cases and provide insights
into their properties. For example, as the complex plane can be tessellated by
regular hexagons and $\overline{\Omega}_{3,2}$ is a regular hexagon, it is
easy to show that $S_{3,2}$ is extendable to an elliptic function; this is
known in \cite{Dixon} and $S_{3,2}$ is the Dixon's elliptic function
$\operatorname*{sm}$. Similarly, it can be shown, based upon Proposition
\ref{(n/2)th Power} and the fact that $\overline{\Omega}_{4,3}$ is a square,
that $(S_{4,3})^{2}$ can be continued to an elliptic function; this is
established in \cite{Levin} via an entirely different (and less elementary) approach.

\item The analyses given here can be equally applied even if $k$ is
\textit{any} positive real number less than $n$. Thus, for any real $p>1$ and
any integral $n\geq2$, we have essentially found the analytic continuation of
the function $y\mapsto\int_{0}^{y}(1-x^{n})^{-1/p}dx$ and of its inverse. We
leave the details to the interested reader.

\item Treating $S_{n,k}$ as a complex function allows us to determine the
maximal interval on which the real function $S_{n,k}$ is real-analytic. The
details can be found in \cite{Ding}.
\end{enumerate}

\end{document}